\documentclass{amsart}
\usepackage{amsfonts,amssymb,epsfig}
\usepackage[latin1]{inputenc}
\newtheorem{thm}{Theorem}[section]

\theoremstyle{definition}
\newtheorem{defn}[thm]{Definition}
\theoremstyle{remark}

\numberwithin{equation}{section}

\newcommand{\cF}{\mathcal F}
\newcommand{\cX}{\mathcal X}

\newcommand{\bbR}{\mathbb R}
\newcommand{\bbT}{\mathbb T}

\newcommand{\ddim}{{\rm ddim\ }}
\newcommand{\rank}{{\rm rank\ }}

\begin{document}

\title{Reduction and Integrability}

\author{Nguyen Tien Zung}
\address{GTA, UMR 5030 CNRS, Département de Mathématiques, Université Montpellier II}
\email{tienzung@math.univ-montp2.fr {\it URL}: www.math.univ-montp2.fr/\~{}tienzung}

\date{first version, 11/Jan/2002}%

\begin{abstract}
We discuss the relationship between the integrability of a dynamical system
invariant under a Lie group action and its reduced integrability, i.e. integrability
of the corresponding reduced system
\end{abstract}
\maketitle

\section{Introduction}

The aim of this note is to address the following question (Q):

\begin{quote}
Given a manifold $M$ (eventually with singularities), a Lie group $G$ which acts on
$M$ (in such a way that the quotient $M/G$ is a manifold with singularities), and a
vector field $X$ on $M$ which is invariant under the action of $G$. Denote by $X/G$
the projection of $X$ on $M/G$. What is the relationship between the integrability
of $(M,X)$ and the integrability of $(M/G, X/G)$ (a.k.a. the {\it reduced
integrability} of $(M,X,G)$) ?
\end{quote}

The above question is very natural, since dynamical systems often admit natural
symmetry groups, and by integrability of a problem in classical mechanics one often
means its reduced integrability. It seems to me, however, that Question (Q) has not
been formally addressed anywhere in the literature, and that's why this note.

We will consider two different cases: Hamiltonian and non-Hamiltonian. For
simplicity, we will assume that $G$ is a compact group. We will show that, when the
action of $G$ is Hamiltonian, Hamiltonian integrability is the same as reduced
Hamiltonian integrability. In the non-Hamiltonian case, integrability still implies
reduced integrability, though the inverse needs not be true. The proof of these
facts is elementary : we simply play with the dimensions of various spaces and their
intersections, quotients, etc.

\section{Hamiltonian integrability}

Let $(M,\Pi)$ be a real Poisson manifold, with $\Pi$ being the Poisson structure.
Let $H$ be a (smooth or analytic) function on $M$, and $X_H$ be the corresponding
Hamiltonian vector field. Assume that we have found a set $\cF$ of first integrals
of $X_H$, i.e. each $F \in \cF$ is a function on $M$ which is preserved by $X_H$
(equivalently, $\{F, H\} = 0$, where $\{.,.\}$ denotes the Poisson bracket as
usual). Denote by $\ddim \cF$ the functional dimension of $\cF$, i.e. the maximal
number of functionally independent functions in $\cF$. To avoid pathologies, we will
always assume that the functional dimension of the restriction of $\cF$ to any open
subset of $M$ is equal to $\ddim \cF$. (This is automatic in the analytic case).

We will associate to $\cF$ the space $\cX = \cX_\cF$ of Hamiltonian vector fields
$X_F$ such that $X_F(G) = 0$ for all $G \in \cF$ and $F$ is functionally dependent
of $\cF$ (i.e. the functional dimension of the union of $\cF$ with the function $F$
is the same as the functional dimension of $\cF$). Clearly, $X_H$ belongs to $\cX$,
and the vector fields in $\cX$ commute pairwise. Denote by $\ddim \cX$ the
functional dimension of $\cX$, i.e. the maximal number of vector fields in $X$ whose
exterior (wedge) product does not vanish.

With the above notations, we have the following definition, due essentially to
Nekhoroshev \cite{Nekhoroshev} and Mischenko and Fomenko \cite{MiFo} :

\begin{defn}
\label{defn:HI}
A Hamiltonian vector field $X_H$ on an $m$-dimensional Poisson
manifold $(M,\Pi)$ is called {\it Hamiltonianly  integrable} with the aid of a set
of first integrals $\cF$, if $m = \ddim \cF + \ddim \cX_\cF$. It is called {\it
properly Hamiltonianly integrable} if $\cF$ satisfies the following additional
properness condition :

There are $q$ functions $F_1,...,F_q$ in $\cF$, where $q = \ddim \cF$, which are
functionally independent, and whose joint ``moment map'' $(F_1,...,F_q) : M \to
\bbR^q$ (here we consider the real case) is a proper map from $M$ to its image, and
there are $p$ vector fields $X_1,...,X_p$ in $\cX$, where $p = \ddim \cX$, such that
the image of their singular set $\{x \in M , X_1 \wedge X_2 \wedge ... \wedge X_p
(x) = 0\}$ under the map $(F_1,...,F_q) : M \to \bbR^q$ is nowhere dense in
$\bbR^q$.
\end{defn}

{\it Remarks}

1. The above notion of integrability is often called {\it generalized Liouville
integrability}, or also {\it non-commutative integrability} by Mischenko-Fomenko,
due to the fact that the functions in $\cF$ do not Poisson-commute in general, and
in many cases one may choose $\cF$ to be a finite-dimensional non-commutative Lie
algebra of functions (under the Poisson bracket). When the functions in $\cF$
Poisson-commute, we get back to the classical integrability à la Liouville.

2. We always have $m \leq \ddim \cF + \ddim \cX$ (even for non-integrable systems),
because the vector fields in $\cX$ are tangent to the common level sets of the
functions in $\cF$. Thus the integrability condition $m = \ddim \cF + \ddim \cX$ is
a maximality, or fullness, condition on $\cF$.

3. Casimir functions of $(M,\Pi)$, i.e. functions whose Hamiltonian vector fields
vanish, must be functionally dependent of $\cF$ in the integrable case - otherwise
we could add them to $\cF$ to increase the functional dimension of $\cF$, which
contradicts the above remark about $m \leq \ddim \cF + \ddim \cX$.

4. It follows directly from the fact that if $X_F \in \cX$ then $F$ is functionally
dependent of $\cF$ that we have $\ddim \cX \leq \ddim \cF$. We didn't mention the
rank of the Poisson structure $\Pi$ in the above definition, but of course $\ddim
\cX \leq 1/2\ \rank \Pi$. When $m = \ddim \cF + \ddim \cX$ and $\ddim \cX < 1/2\
\rank \Pi$, some authors also say that the system is {\it super-integrable} (because
in this case one finds more first integrals than necessary for integrability).

5. We prefer to use the term {\it Hamiltonian integrability}, to contrast it with
the {\it non-Hamiltonian integrability} discussed in the next section.

6. Under the additional properness condition, one get a natural generalization of
the classical Liouville theorem \cite{Nekhoroshev,MiFo}: the manifold $M$ is
foliated by invariant isotropic tori on which the flow of $X_H$ is quasi-periodic
(thus the behavior of $X_H$ is very regular, justifying the word ``integrable''),
and there also exist local generalized action-angle coordinates. The existence of
action-angle coordinates for Liouville-integrable systems is often referred to as
Arnold-Liouville theorem, though it was probably first proved by Mineur
\cite{Mineur}.

Assume that there is a Hamiltonian action of a Lie group $G$ on $(M,\Pi)$, given by
an equivariant moment map $\pi : M \to {\frak g}^{\ast}$, where ${\frak g}$ denotes
the Lie algebra of $G$, such that the following conditions are satisfied :

1) The action of $G$ on $M$ is proper, so that the quotient space $M/G$ is a
singular manifold whose ring of functions may be identified with the ring of
$G$-invariant functions on $M$.

2) Recall that the image $\pi(M)$ of $M$ under the moment map $\pi : M \to {\frak
g}^{\ast}$ is saturated by symplectic leaves (i.e. coadjoint orbits) of ${\frak
g}^{\ast}$. Denote by $s$ the minimal codimension in ${\frak g}^{\ast}$ of a
coadjoint orbit which lies in $\pi(M)$. Then we assume that there exist $s$
functions $f_1,...,f_s$ on ${\frak g}^{\ast}$, which are invariant on the coadjoint
orbits which lie in $\pi(M)$, and such that for almost every point $x \in M$ we have
$df_1 \wedge ... \wedge df_s (\pi(x)) \neq 0$.

For example, when $G$ is compact and $M$ is connected, then the above conditions are
satisfied automatically.

If $\pi(M)$ contains a generic point of ${\frak g}^{\ast}$, then $s = {\rm ind\ }
{\frak g}$, where ${\rm ind\ } {\frak g}$ denotes the index of $\frak g$, i.e. the
corank of the corresponding linear Poisson structure on ${\frak g}^{\ast}$, and if
there are $({\rm ind\ } {\frak g})$ global functionally independent Casimir
functions on ${\frak g}^{\ast}$ they may play the role of required functions
$f_1,...,f_s$ in the above assumption. If $\pi(M)$ lies in the singular part of
${\frak g}^{\ast}$ then $s$ may be greater than the index of $\frak g$.

Since $\Pi$ is preserved by $G$, it can be projected to a Poisson structure, denoted
by $\Pi/G$, on $M/G$. In the case when $M$ is a symplectic manifold, the symplectic
leaves of $(M/G,\Pi/G)$ are known as Marsden-Weinstein reductions.

Let $H$ be a function on $M$ which is invariant by $G$. Then $H$ may be viewed as
the pull-back of a function $h$ on $M/G$ via the projection ${\frak p}: M \to M/G$,
$H = {\frak p}^{\ast}(h)$. Denote by $X_H$ (resp., $X_h$) the Hamiltonian vector
field of $H$ (resp., $h$) on $(M,\Pi)$ (resp., $(M/G,\Pi/G)$). Of course, $X_H$ is
$G$-invariant, and its projection to $M/G$ is $X_h$.

With the above notations and assumptions, we have :

\begin{thm}
\label{thm:HI1} If the system $(M/G,X_h)$ is Hamiltonianly integrable, then the
system $(M,X_H)$ also is Hamiltonianly integrable. If $G$ is compact and $(M/G,X_h)$
is properly Hamiltonianly integrable, then $(M,X_H)$ also is properly Hamiltonianly
integrable.
\end{thm}

{\it Proof}. Denote by $\cF'$ a set of first integrals of $X_h$ on $M/G$ which
provides the integrability of $X_h$, and by $\cX' = \cX_{\cF'}$ the corresponding
space of commuting Hamiltonian vector fields on $M/G$. We have $\dim M/G = p' + q'$
where $p' = \ddim \cX'$ and $q' = \ddim \cF'$.

Recall that, by our assumptions, there exist $s$ functions $f_1,...,f_s$ on ${\frak
g}^{\ast}$, which are functionally independent almost everywhere in $\pi(M)$, and
which are invariant on the coadjoint orbits which lie in $\pi(M)$. Here $s$ is the
minimal codimension in ${\frak g}^{\ast}$ of the coadjoint orbits which lie in
$\pi(M)$. We can complete $(f_1,...,f_s)$ to a set of $d$ functions
$f_1,...,f_s,f_{s+1},...,f_d$ on ${\frak g}^{\ast}$, where $d = \dim G = \dim \frak
g$ denotes the dimension of $\frak g$ , which are functionally independent almost
everywhere in $\pi(M)$.

Denote by $\overline{\cF}$ the pull-back of $\cF'$ under the projection ${\frak p} :
M \to M/G $, and by $F_1,...,F_d$ the pull-back of $f_1,...,f_d$ under the moment
map $\pi : M \to {\frak g}^{\ast}$. Note that, since $H$ is $G$-invariant, the
functions $F_i$ are first integrals of $X_H$. And of course, $\overline{\cF}$ is
also a set of first integrals of $X_H$. Denote by $\cF$ the union of
$\overline{\cF}$ with $(F_{s+1},...,F_d)$. (It is not necessary to include
$F_1,...,F_s$ in this union, because these functions are $G$-invariant and project
to Casimir functions on $M/G$, which implies that they are functionally dependent of
$\overline{\cF}$). We will show that $X_H$ is Hamiltonianly integrable with the aid
of $\cF$.

Notice that, by assumptions, the coadjoint orbits of ${\frak g}^{\ast}$ which lie in
$\pi(M)$ are of generic dimension $d-s$, and the functions $f_{s+1},...,f_d$ may be
viewed as a coordinate system on a symplectic leaf of $\pi(M)$ at a generic point.
In particular, we have
$$< df_{s+1} \wedge ... \wedge df_d, X_{f_{s+1}} \wedge ...
X_{f_d} > \neq 0,$$ which implies, by equivariance :
$$< dF_{s+1} \wedge ... \wedge dF_d, X_{F_{s+1}} \wedge ... X_{F_d} > \neq 0.$$

Since the vector fields $X_{F_{s+1}},...,X_{F_d}$ are tangent to the orbits of $G$
on $M$, and the functions in $\overline{\cF}$ are invariant on the orbits of $G$, it
implies that the set $(F_{s+1},...,F_{d})$ is ``totally'' functionally independent
of $\overline{\cF}$. In particular, we have :

\begin{equation}
\ddim \cF = \ddim \cF' + \ddim (F_{s+1},...,F_{d}) = q' + d - s ,
\end{equation}
where $q' = \ddim \cF'$. On the other hand, we have

$$\dim M = \dim M/G + (d-k) = p' + q' + d - k ,$$
where $p' = \ddim \cX_{\cF'}$, and $k$ is the dimension of a minimal isotropic group
of the action of $G$ on $M$. Thus, in order to show the integrability condition
$$\dim M = \ddim \cF + \ddim \cX_{\cF} ,$$ it remains to show that

\begin{equation}
\label{eqn:ddimcX}
\ddim \cX_{\cF} = \ddim \cX_{\cF'} + (s-k) .
\end{equation}

Consider the vector fields $Y_1 = X_{F_1},...,Y_d = X_{F_d}$ on $M$. They span the
tangent space to the orbit of $G$ on $M$ at a generic point. The dimension of such a
generic tangent space is $d-k$. It implies that, among the first $s$ vector fields,
there are at least $s-k$ vector fields which are linearly independent at a generic
points : we may assume that $Y_1 \wedge ... \wedge Y_{s-k} \neq 0 .$

Let $X_{h_1},...,X_{h_{p'}}$ be $p'$ linearly independent (at a generic point)
vector fields which belong to $\cX_{\cF'}$, where $p' = \ddim \cX_{\cF'}$. Then we
have

$$
X_{{\frak p}^{\ast}(h_1)},...,X_{{\frak p}^{\ast}(h_{p'})},Y_1,...,Y_{s-k} \in
\cX_{\cF},
$$
and these $p' + s - k$ vector fields are linearly independent at a generic point.
(Recall that, at each point $x \in M$, the vectors $Y_1(x),...,Y_{s-k}(x)$ are
tangent to the orbit of $G$ which contains $x$, while the linear space spanned by
$X_{{\frak p}^{\ast}(h_1)},...,X_{{\frak p}^{\ast}(h_{p'})}$ contains no tangent
direction to this orbit).

Thus we have $\ddim \cX_{\cF} \geq p' + s - k$, which means that $\ddim \cX_{\cF} =
p' + s - k$ (because, as discussed earlier, we always have $\ddim \cF + \ddim
\cX_\cF \leq \dim M$). We have proved that if $(M/G,X_h)$ is Hamiltonianly
integrable then $(M,X_H)$ also is.

Now assume that $G$ is compact and $(M/G,X_h)$ is properly Hamiltonianly integrable
: there are $q'$ functionally independent functions $g_1,...,g_{q'} \in \cF'$ such
that $(g_1,...,g_{q'}): M/G \to \bbR^{q'}$ is a proper map from $M/G$ to its image,
and $p'$ Hamiltonian vector fields $X_{h_1},...,X_{h_{p'}}$ in $\cX'$ such that on a
generic common level set of $(g_1,...,g_{q'})$ we have that $X_{h_1} \wedge ...
\wedge X_{h_{p'}}$ does not vanish anywhere. Then it is straightforward that
$${\frak p}^{\ast}(g_1),...,{\frak p}^{\ast}(g_{q'}),F_{s+1},...,F_d \in \cF$$ and
the map $$({\frak p}^{\ast}(g_1),...,{\frak p}^{\ast}(g_{q'}),F_{s+1},...,F_d) : M
\to \bbR^{q'+d-s}$$ is a proper map from $M$ to its image. More importantly, on a
generic level set of this map we have that the $(q'+s-k)$-vector $X_{{\frak
p}^{\ast}(h_1)} \wedge ... \wedge X_{{\frak p}^{\ast}(h_{p'})} \wedge Y_1 \wedge ...
\wedge Y_{s-k}$ does not vanish anywhere. To prove this last fact, notice that
$X_{{\frak p}^{\ast}(h_1)} \wedge ... \wedge X_{{\frak p}^{\ast}(h_{p'})} \wedge Y_1
\wedge ... \wedge Y_{s-k} (x) \neq 0$ for a point $x \in M$ if and only if
$X_{{\frak p}^{\ast}(h_1)} \wedge ... \wedge X_{{\frak p}^{\ast}(h_{p'})} (x) \neq
0$ and $Y_1 \wedge ... \wedge Y_{s-k} (x) \neq 0$ (one of these two multi-vectors is
transversal to the $G$-orbit of $x$ while the other one ``lies on it''), and that
these inequalities are $G \times \bbR^{p'}$-invariant properties, where the action
of $\bbR^{p'}$ is generated by $X_{{\frak p}^{\ast}(h_1)},...,X_{{\frak
p}^{\ast}(h_{p'})}$. $\diamondsuit$ \\

{\it Remark}. Recall from Equation (\ref{eqn:ddimcX}) above that we have $\ddim
\cX_{\cF} - \ddim \cX_{\cF'} = s-k$, where $k$ is the dimension of a generic
isotropic group of the $G$-action on $M$, and $s$ is the (minimal) corank in ${\frak
g}^{\ast}$ of a coadjoint orbit which lies in $\pi(M)$. On the other hand, the
difference between the rank of the Poisson structure on $M$ and the reduced Poisson
structure on $M/G$ can be calculated as follows :
\begin{equation}
\rank \Pi - \rank \Pi/G = (d-k) + (s-k)
\end{equation}
Here $(d-k)$ is the difference between $\dim M$ and $\dim M/G$, and $(s-k)$ is the
difference between the corank of $\Pi/G$ in $M/G$ and the corank of $\Pi$ in $M$. It
follows that
\begin{equation}
\rank \Pi - 2 \ddim \cX_{\cF} = \rank \Pi/G - 2 \ddim \cX_{\cF'} + (d-s)
\end{equation}
In particular, if $d-s > 0$ (typical situation when $G$ is non-Abelian), then we
always have $\rank \Pi - 2 \ddim \cX_{\cF} > 0$ (because we always have $\rank \Pi/G
- 2 \ddim \cX_{\cF'} \geq 0$ due to integrability), i.e. the original system is
always super-integrable with the aid of $\cF$. When $G$ is Abelian (implying $d=s$),
and the reduced system is Liouville-integrable with the aid of $\cF'$ (i.e. $\rank
\Pi/G = 2 \ddim \cX_{\cF'}$), then the original system is also Liouville-integrable
with the aid of $\cF$.

{\it Remark.} Following Mischenko-Fomenko \cite{MiFo}, we will say that a
hamiltonian system $(M,\Pi,X_H)$ is {\it non-commutatively integrable in the
restricted sense} with the aid of $\cF$ , if $\cF$ is a finite-dimensional Lie
algebra under the Poisson bracket and $(M,\Pi,X_H)$ is Hamiltonianly integrable with
the aid of $\cF$. In other words, we have an equivariant moment maps $(M,\Pi) \to
{\frak f}^{\ast}$, where ${\frak f}$ is some finite-dimensional Lie algebra, and if
we denote by $f_1,...,f_n$ the components of this moment map, then they are first
integrals of $X_H$, and $X_H$ is Hamiltonianly integrable with the aid of this set
of first integrals. Theorem \ref{thm:HI1} remains true, and its proof remains the
same if not easier, if we replace Hamiltonian integrability by non-commutative
integrability in the restricted sense. Indeed, if $M \to {\frak g}^{\ast}$ is the
equivariant moment map of the symmetry group $G$, and if $M/G \to {\frak h}^{\ast}$
is an equivariant moment map which provides non-commutative integrability in the
restricted sense on $M/G$, then the map $M \to {\frak h}^{\ast}$ (which is the
composition $M \to M/G \to{\frak h}^{\ast}$) is an equivariant moment map which
commutes with $M \to {\frak g}^{\ast}$, and the direct sum of this two maps, $M \to
{\frak f}^{\ast}$ where $\frak f = \frak g \bigoplus \frak h$, will provide
non-commutative integrability in
the restricted sense on $M$. \\

The above remarks show that {\it the notion of Hamiltonian integrability (or
non-commutative integrability if you prefer), rather than integrability à la
Liouville, is the most natural one when dealing with systems admitting (non-Abelian)
symmetry groups}. \\

{\it Examples}. 1) The simplest example which shows an evident relationship between
reduction and integrability is the classical Euler top : it can be written as a
Hamiltonian system on $T^{\ast}SO(3)$, invariant under a natural Hamiltonian action
of $SO(3)$, is integrable with the aid of a set of four first integrals, and with
2-dimensional isotropic invariant tori. 2) The geodesic flow of a bi-invariant
metric on a compact Lie group is also properly Hamiltonianly integrable :
in fact, the corresponding reduced system is trivial (identically zero). \\

It is known that most Hamiltonianly integrable systems are also integrable à la
Liouville, i.e. there exists a set of first integrals which commute pairwise and
which make the system integrable - see e.g. \cite{Fomenko} for a very long
discussion on this subject. A question of similar kind, which is directly related to
the inverse of Theorem \ref{thm:HI1}, is the following :

\begin{quote}
If $X_H$ is Hamiltonianly integrable with the aid of $\cF_1$, and if $\cF_2$ is
another set of first integrals of $X_H$ which contains $\cF_1$, then is it true that
$X_H$ is also Hamiltonianly integrable with the aid of $\cF_2$ ? In particular, let
$\cF_H$ denotes the set of all first integrals of $X_H$. If $X_H$ is Hamiltonianly
integrable, then is it true that it is integrable with the aid of $\cF_H$ ?
\end{quote}

A related question is the following :

\begin{quote}
Suppose that $\cF = (F_1,...,F_q)$ is a set of independent first integrals of $X_H$
such that regular common level sets of these functions $(F_1,...,F_q)$ are isotropic
submanifolds of $(M,\Pi)$. Is it true that $X_H$ is integrable with the aid of $\cF$
?
\end{quote}

Remark that the inverse to the later question is always true. And if we can say YES
to the later question, then we can also say YES to the former one, because adding
first integrals has the affect of minimizing invariant submanifolds, and a
submanifold of an isotropic submanifold is again an isotropic submanifold. It is
easy to see that, at least in the smooth proper case, the answer to the above two
questions is YES : smooth proper Hamiltonian integrability of a Hamiltonian system
is equivalent to the singular foliation of the Poisson manifold by invariant
isotropic tori (This fact is similar to the classical Liouville theorem).

For the following theorem, we use the same notations and preliminary assumptions as
in Theorem \ref{thm:HI1} :

\begin{thm}
\label{thm:HI2} If $G$ is compact, and if the Hamiltonian system $(M,X_H)$ is
Hamiltonianly integrable with the aid of $\cF_H$ (the set of all first integrals),
then the reduced Hamiltonian system $(M/G,X_h)$ is also Hamiltonianly integrable.
The same thing holds in the smooth proper case.
\end{thm}

{\it Proof}. By assumptions, we have $\dim M = p + q$, where $q = \ddim \cF_H$ and
$p = \ddim \cX_{\cF_H}$, and we can find $p$ first integrals $H_1,...,H_p$ of $H$
such that $X_{H_1},...,X_{H_p}$ are linearly independent (at a generic point) and
belong to $\cX_{\cF_H}$. In particular, we have $X_{H_i}(F) = 0$ for any  $F \in
\cF$ and $1 \leq i \leq p$.

An important observation is that the functions $H_1,...,H_p$ are $G$-invariant. In
deed, if we denote by $F_1,...,F_d$ the components of the equivariant moment map
$\pi : M \to {\frak g}^{\ast}$ (via an identification of ${\frak g}^{\ast}$ with
$\bbR^d$), then since $H$ is $G$-invariant we have $\{H,F_j\}=0$, i.e. $F_j \in
\cF_H$, which implies that $\{F_j,H_i\} = 0 \ \ \forall 1 \leq i \leq d, \ \ 1 \leq
j \leq p$, which means that $H_i$ are $G$-invariant.

Denote by $h_i$ the projection of $H_i$ on $M/G$ (recall that the projection of $H$
on $M/G$ is denoted by $h$). Then the Hamiltonian vector fields $X_{h_i}$ belong to
$\cX_{\cF_h}$ : Indeed, if $f \in \cF_h$ then ${\frak p}^{\ast}(f)$ is a first
integral of $H$, implying $\{H_i,{\frak p}^{\ast}(f)\} = 0$, or $\{h_i,f\}=0$, where
$\frak p$ denotes the projection $M \to M/G$.

To prove the integrability of $X_h$, it is sufficient to show that
\begin{equation}
\label{eqn:Xh} \dim M/G \leq \ddim \cF_h + \ddim (X_{h_1},...,X_{h_q})
\end{equation}
But we denote by $r$ the generic dimension of the intersection of a common level set
of $p$ independent first integrals of $X_H$ with an orbit of $G$ in $M$, then one
can check that
$$
p - \ddim (X_{h_1},...,X_{h_q}) = \ddim \cX_{\cF_H} - \ddim (X_{h_1},...,X_{h_q}) =
r
$$
and
$$ q - \ddim \cF_h = \ddim \cF_H - \ddim \cF_h \leq (d-k) - r
$$
where $(d-k)$ is the dimension of a generic orbit of $G$ in $M$. To prove the last
inequality, notice that functions in $\cF_h$ can be obtained from functions in
$\cF_H$ by averaging with respect to the $G$-action. Also, $G$ acts on the
(separated) space of common level sets of the functions in $\cF_H$, and isotropic
groups of this $G$-action are of (generic) codimension $(d-k) - r$.

The above two formulas, together with $p+q = \dim M = \dim M/G + (d-k)$, implies
Inequality (\ref{eqn:Xh}) (it is in fact an equality).

We will leave the proper case to the reader as an exercise. $\diamondsuit$

\section{Non-Hamiltonian integrability}

The interest in non-Hamiltonian integrability comes partly from the fact that there
are many non-Hamiltonian (e.g. non-holonomic) systems whose behaviors are very
similar to that of integrable Hamiltonian systems, see e.g. \cite{BaCu,CuDu}. In
particular, in the proper case, the manifold is foliated by invariant tori on each
of which the system is quasi-periodic. Another common point between (integrable)
Hamiltonian and non-Hamiltonian systems is that their local normal form theories are
very similar and are related to local torus actions, see e.g.
\cite{Stolovitch,ZungNormal}. The notion of non-Hamiltonian integrability was
probably first introduced by Bogoyavlenskij \cite{Bogoyavlenskij}, who calls it {\it
broad integrability}, in his study of tensor invariants of dynamical systems. Let us
give here a definition of non-Hamiltonian integrability, which is similar to the
ones found in \cite{BaCu,Bogoyavlenskij,CuDu,Stolovitch,ZungNormal} :

\begin{defn}
A vector field on a (eventually singular) manifold $M$ is called {\it
non-Hamiltonianly integrable} with the aid of $(\cF,\cX)$, where $\cF$ is a set of
funtions on $M$ and
$\cX$ is a set of vector fields on $M$, if the following conditions are satisfied : \\
a) Functions in $\cF$ are first integrals of $X$: $X(F) = 0 \ \ \forall F \in \cF .$
\\
b) Vector fields in $\cX$ commute pairwise and commute with $X$ : $[Y,Z] = [Y,X] = 0
\ \ \forall \ Y,Z \in \cX .$ \\
c) Functions in $\cF$ are common first integrals of vector fields in $X$: $Y(F) = 0
\ \ \forall \ Y \in \cX, F \in \cF .$ \\
d) $\dim M = \ddim \cF + \ddim \cX .$ \\
If, moreover, there exist $p$ vector fields $Y_1,...,Y_p \in \cX$ and $q$
functionally independent functions $F_1,...,F_q \in \cF$, where $p = \ddim \cX$ and
$q = \ddim \cF$, such that the map $(F_1,...,F_q) : M \to \bbR^q$ is a proper map
from $M$ to its image, and for almost any level set of this map the vector fields
$Y_1,...,Y_p$ are linearly independent everywhere on the level set, then we say that
$X$ is {\it properly non-Hamiltonianly integrable} with the aid of $(\cF,\cX)$.
\end{defn}

{\it Remarks.}

1. It is straightforward that, in the proper case,  the manifold is a (singular)
foliation by invariant tori (common level sets of some first integrals) on each of
which the vector field $X$ is quasi-periodic. (This is similar to the classical
Liouville theorem).

2. If a Hamiltonian system is (properly) Hamiltonianly integrable, then it is also
(properly) non-Hamiltonianly integrable, though the inverse is not true : it may
happen that the invariant tori are not isotropic, see e.g.
\cite{Bogoyavlenskij,Fasso} for a detailed discussion about this question. \\

One of the main differences between the non-Hamiltonian case and the Hamiltonian
case is that reduced non-Hamiltonian integrability does not imply integrability. In
fact, in the Hamiltonian case, we can lift Hamiltonian vector fields from $M/G$ to
$M$ via the lifting of corresponding functions. In the non-Hamiltonian case, no such
canonical lifting exists, therefore commuting vector fields on $M/G$ do not provide
commuting vector fields on $M$. For example, consider a vector field of the type $X
= a_1 \partial/\partial x_1 + a_2 \partial/\partial x_2 + b(x_1,x_2)
\partial/\partial x_3$ on the standard torus $\bbT^3$ with periodic coordinates
$(x_1,x_2,x_3)$, where $a_1$ and $a_2$ are two incommensurable real numbers
($a_1/a_2 \notin {\mathbb Q}$), and $b(x_1,x_2)$ is a smooth function of two
variables. Then clearly $X$ is invariant under the ${\mathbb S}^1$-action generated
by $\partial/\partial x_3$, and the reduced system is integrable. On the other hand,
for $X$ to be integrable, we must be able to find a function $c(x_1,x_2)$ such that
$[X, \partial/\partial x_1 + c(x_1,x_2)\partial/\partial x_3] = 0 $. This last
equation does not always have a solution (it is a small divisor problem, and depends
on $a_1/a_2$ and the behavior of the coefficients of $b(x_1,x_2)$ in its Fourier
expansion), i.e. there are choices of $a_1,a_2,b(x_1,x_2)$ for which the vector
field $X$ is not integrable.
\\

However, non-Hamiltonian integrability still implies reduced integrability. Before
formulating a precise result, let us mention a question similar to the one already
mentioned in the previous section :

\begin{quote}
For a vector field $X$ on a manifold $M$, denote by $\cF_X$ the set of all first
integrals of $X$, and by $\cX_X$ the set of vector fields which preserve each
function in $\cF$ and commute with $X$. Suppose that $X$ is non-Hamiltonianly
integrable. Is it then non-Hamiltonianly integrable with the aid of $(\cF_X,\cX_X)$
? In other words, is it true that vector fields in $\cF$ commute pairwise and $\ddim
\cX_X + \ddim \cF_X = \dim M$ ?
\end{quote}

It is easy to see that the answer to the above question is YES in the proper
non-Hamiltonianly integrable case, under the additional assumption that the orbits
of $X$ are dense (i.e. its frequencies are incommensurable) on almost every
invariant torus (i.e. common level of a given set of first integrals $\cF$). In this
case $\cX_X$ consists of the vector fields which are quasi-periodic on each
invariant torus. Another case where the answer is also YES arises in the study of
local normal forms of analytic integrable vector fields, see e.g. \cite{ZungNormal}.

\begin{thm}
\label{thm:nHI} Let $X$ be a smooth properly non-Hamiltonianly integrable system on
a manifold $M$ with the aid of $(\cF_X,\cX_X)$, and $G$ be a compact Lie group
acting on $M$ which preserves $X$. Then the reduced system on $M/G$ is also properly
non-Hamiltnianly integrable.
\end{thm}

{\it Proof.} Let $\cX_X^G$ denote the set of vector fields which belong to $\cX_X$
and which are invariant under the action of $G$. Note that the elements of $\cX_X^G$
can be obtained from the elements of $\cX_X$ by averaging with respect to the
$G$-action.

A key ingredient of the proof is the fact $\ddim \cX_X^G = \ddim \cX_X$ (To see this
fact, notice that near each regular invariant torus of the system there is an
effective torus action (of the same dimension) which preserves the system, and this
torus action must necessarily commute with the action of $G$. The generators of this
torus action are linearly independent vector fields which belong to  $\cX_X^G$ - in
fact, they are defined locally near the union of $G$-orbits which by an invariant
torus, but then we can extend them to global vector fields which lie in $\cX_X^G$)

Therefore, we can project the pairwise commuting vector fields in $\cX_X^G$ from $M$
to $M/G$ to get pairwise commuting vector fields on $M/G$. To get the first
integrals for the reduced system, we can also take the first integrals of $X$ on $M$
and average them with respect to the $G$-action to make them $G$-invariant.  The
rest of the proof of Theorem
\ref{thm:nHI} is similar to that of Theorem \ref{thm:HI2}. $\diamondsuit$ \\


{\bf Acknowledgements}. A part of this note, Theorem \ref{thm:HI1}, dates back to
1999-2000 when I was preparing some lectures on the subject and was surprised by the
lack of such a theorem in the standard literature. I would like to thank A.T.
Fomenko for encouraging me to write up this note.

\bibliographystyle{amsplain}

\end{document}